\documentclass{article}

\usepackage[top=1.5in, bottom=1.5in, left=1.7in, right=1.7in]{geometry}

\title{Element order versus minimal degree in permutation groups: 
         an old lemma with new applications}
\author{L\'aszl\'o Babai\thanks{University of Chicago.
 \href{http://people.cs.uchicago.edu/\~laci/}{people.cs.uchicago.edu/$\sim$laci}}
\ and \'Akos Seress\thanks{Ohio State University.  Deceased 
 February 13, 2013.}
}

\usepackage{amsmath}
\usepackage{amsthm}
\usepackage{amssymb}
\usepackage{hyperref}

\theoremstyle{definition} \newtheorem{definition}{Definition}[section]
\theoremstyle{definition} \newtheorem{conjecture}{Conjecture}[section]
\theoremstyle{plain} \newtheorem{theorem}[definition]{Theorem}
\theoremstyle{plain} 
\theoremstyle{plain} \newtheorem{proposition}[definition]{Proposition}
\theoremstyle{plain} 
\theoremstyle{plain} 
\DeclareMathOperator{\supp}{supp}

\begin{document}


\maketitle

\begin{abstract}
In this note we present a simplified and slightly generalized version
of a lemma the authors published in 1987.  The lemma as stated here
asserts that if the order of a permutation of $n$ elements is
greater than $n^{\alpha}$ then some nonidentity power of the
permutation has support size less than $n/\alpha$.  The original
version made an unnecessary additional assumption on the cycle
structure of the permutation; the proof of the present cleaner
version follows the original proof verbatim.   Application areas 
include parallel and sequential algorithms for permutation groups,
the diameter of Cayley graphs of permutation groups, and the automorphisms of
combinatorial structures with regularity constraints such as 
Latin squares, Steiner 2-designs, and strongly regular graphs.
This note also serves as a modest tribute to
the junior author whose untimely passing is deeply mourned.
\end{abstract}

\section{The lemma}

\noindent
For a permutation $\pi$ on the set $K$, the \emph{support}
of $\pi$ is defined as the set $\supp(\pi) = \{x\in K \mid x^{\pi}\neq x\}$.
$S_n$ denotes the symmetric group of degree $n$ (and order $n!$).
The \emph{degree} of a permutation is the size of its support, and the
\emph{minimal degree} of a permutation group is the smallest degree of
its non-identity elements.  The following lemma shows that if there are
elements of large order in a permutation group then its minimal degree 
is small.

\medskip\noindent
{\bf Lemma.}
Let $\pi\in S_n$ have order $\ge n^{\alpha}$ for some $\alpha > 0$.
Then $|\supp(\pi^m)|\le n/\alpha$ for some power $\pi^m\neq 1$.

\medskip\noindent
\begin{proof}
Let $\pi$ act on the set $K$ where $|K|=n$.  
Let the order of $\pi$ be $N = n^{\alpha} =\prod_{i=1}^r q_i$
where $q_i= p_i^{\beta_i} >1$ are powers of distinct primes $p_i$.
For each $x\in K$, let us consider the set $P(x)$ of those $i$
for which $q_i$ divides the length of the $\pi$-cycle through $x$.
Clearly, for each $x\in K$,
\begin{equation}   \label{eq:cycle}
       \prod_{i\in P(x)} q_i \le n .
\end{equation}
Let $n(i)$ denote the number of points $x\in K$ such that 
$i\in P(x)$.  Let us estimate the weighted average $W$ of the $n(i)$
with weights $\log q_i$.  Recall that the sum of weights is 
$\sum \log q_i = \log N \ge \alpha \log n$, therefore
(using Eq.~\eqref{eq:cycle})
\begin{equation}
   W \le \sum_{x\in K} \sum_{i\in P(x)} \log q_i/(\alpha \log n)
     \le (n\log n)/(\alpha\log n) = n/\alpha.
\end{equation}
Thus we infer that $n(i)\le n/\alpha$ for some $i\le r$.
Now let $m=N/p_i$ be the corresponding maximal divisor of $N$.
Clearly $\pi^m$ is not the identity and it fixes all
but $n(i)$ points.  
\end{proof}

\section{History and applications}

The original version of this lemma, published as 
\cite[Lemma 3]{bs-transitivity}, assumed that the permutation 
$\pi$ involved cycles of prime lengths $p_i$ such that 
$\prod p_i \ge n^{\alpha}$.  As we have seen, this assumption is
unnecessary.  The proof of the above cleaner form requires no new 
ideas, however; it is an essentially verbatim copy of the original proof.

Applications of this lemma are manifold, both old and new.

\subsection{Degree of transitivity, diameter, parallel and sequential
complexity of permutation groups}
In \cite{bs-transitivity}, the authors used this lemma to obtain a
short proof of a theorem of Jordan~\cite{jordan1,jordan2}
on the degree of transitivity of
permutation groups.  In another simultaneous paper
\cite{bs2-Sn-diam}, the authors used this lemma to prove an
$\exp(\sqrt{n\ln n}(1+o(1)))$ upper bound on the diameter of all
Cayley graphs of $S_n$; this bound was not improved until very
recently~\cite{helfgott}.  The breakthrough 2013
paper~\cite{helfgott} by Helfgott and Seress that reduced
the diameter bound to quasipolynomial (exponential in a
constant power of $\log n$) makes substantial use
of a slight generalization of the lemma.  In \cite{bs3-permdiam}
the present authors, again using the Lemma, extended the
$\exp(\sqrt{n\ln n}(1+o(1)))$ bound to all permutation groups of
degree $n$ (subgroups of $S_n$), which in this form is tight since
$S_n$ has cyclic subgroups of order $\exp(\sqrt{n\ln n}(1+o(1)))$ 
(product of small prime length cycles).  In \cite{bls-NC}, the
Lemma played a key role in settling the parallel complexity of the
permutation group membership problem; the subsequent papers
\cite{bls-sequential1,bls-sequential2} 
used the Lemma to design and analyze improved 
sequential algorithms for the same problem.

\subsection{Latin squares}
The present note was prompted by a conversation between Babai and
Ian Wanless in July 2013 at a conference celebrating Peter Cameron
at Queen Mary, University of London.  
Wanless mentioned the following recent result of his with
Brendan McKay and Xiande Zhang.  Recall that a \emph{quasigroup} is a set
with a binary operation such that all equations of the form
$ax=b$ and $ya=b$ are uniquely solvable.  In other words,
a quasigroup is a set with a binary operation of which
the multiplication table is a Latin square.

\begin{theorem}[McKay, Wanless, Zhang~\cite{wanless}]
No automorphism of a quasigroup of order $n$ has order greater than $n^2/4.$
\end{theorem}

\medskip\noindent 
We show that a slightly weaker bound, $n^2$, is immediate from
the Lemma.  Indeed, let $\pi$ be a non-identity automorphism of some 
quasigroup $G$ of order $n$. 
Assume the order of $\pi$ is greater than $n^2$.
Then some power $\pi^m\neq 1$ fixes more than half the elements of
$G$.  But the set of fixed points of a non-identity automorphism is a
proper sub-quasigroup and therefore has order $\le n/2$,
a contradiction.  \qed

\medskip\noindent
The McKay--Wanless--Zhang argument is not dissimilar to
ours.   They in fact prove the following stronger result.

\begin{theorem}[McKay, Wanless, Zhang~\cite{wanless}]
No autotopism of a quasigroup of order $n$ has order greater than $n^2/4.$
\end{theorem}
An \emph{autotopism} of a quasigroup $G$ is a 
triple $(\alpha,\beta,\gamma)$ of permutations of
the set $G$ such that for all $g,h\in G$ we have
$\alpha(g)\beta(h)=\gamma(gh)$.  

Wanless pointed out that a quadratic bound, $9n^2$, on the order of
autotopisms of quasigroups of order $n$ also follows quickly from the
Lemma.  Slightly modifying his argument, we infer a bound of $4n^2$
on the order of any autotopism.  (Of course this is still a factor
of 16 worse than their result.)

Indeed, assume $\theta=(\alpha,\beta,\gamma)$ is an autotopism of
order greater than $(2n)^2$.  We can view the autotopisms as acting on
the union $G_1\, \dot\cup\, G_2\, \dot\cup\, G_3$ of three disjoint
copies of $G$.  In fact, the action on $U=G_1\,\dot\cup\, G_2$ is faithful;
let us consider this action.  By the Lemma, some power $\theta^m\neq 1$ 
has more than $n$ fixed points on $U$.  Therefore it has at least
one fixed point in each of $G_1$ and $G_2$.  A result by McKay,
Meynert, and Myrvold~\cite{myrvold} asserts that if a non-identity
autotopism has a fixed point in $G_1$ and a fixed point in
$G_2$ then the number of fixed points in each part is the same and
that number cannot exceed $n/2$, a contradiction.  \qed

We do not expect the quadratic rate of growth in these results to
be optimal.  Let $f(n)$ denote the maximum of the orders of
automorphisms of quasigroups of order $\le n$.
By the above, we have $f(n)=O(n^2)$.

\begin{conjecture}
  \quad\qquad\qquad\qquad  $f(n) = o(n^2) .$
\end{conjecture}

On the other hand, McKay et al.~\cite{wanless} conjecture that 2 is
the best possible exponent.

\begin{conjecture}[\cite{wanless}]
For every $\epsilon > 0$ and for infinitely many values of $n$,
           $$f(n) > n^{2-\epsilon} .$$
\end{conjecture}

The construction of quasigroups with automorphisms of 
nearly quadratic order seems to face significant obstacles;
currently, no superlinear examples are known (cf.~\cite{wanless}).

It is noted in \cite{wanless} that the quadratic upper bound
on the order of automorphisms holds in particular for 
Steiner Triple Systems (STSs) because such systems can be viewed
as quasigroups: if $x,y$ are distinct points of an STS
then define $xy$ as the third point in the unique triple containing
$x,y$; and set $xx=x$.

We observe that the quadratic upper bound for STSs extends to 
all Steiner 2-designs.   A Steiner 2-design (also called
a ``regular linear space'') is an incidence geometry with 
lines of uniform length and exactly one line through every pair 
of points.   If it has $n$ points and each line has $k$ points
then this is a $2-(n,k,1)$-design. 

\begin{proposition}
Let $X$ be a  Steiner 2-design with $n$ points and with lines of length 
$k \ge 3$.  Let $m$ be the maximum order of automorphisms of $X$.  Then
\begin{itemize}
  \item[(a)]  $m < n^2$
  \item[(b)]  $m = O\left(n^{1+\frac{1}{k-2}}\right)$ where
     the implied constant is absolute.
\end{itemize}
\end{proposition}
The proof follows by combining the Lemma with a bound on the
number of fixed points of automorphisms of a Steiner 2-design
by Davies~\cite{davies}.   The details will appear
in~\cite{bab-srg2}.

\subsection{Strongly regular graphs}
In another paper~\cite{bab-srg}, the Lemma is used to
establish strong structural constraints on the automorphism groups 
of strongly regular (SR) graphs.   

Recall that a SR
graph with parameters $(n,k,\lambda,\mu)$ has $n$ vertices,
is regular of degree $k$, each pair of adjacent vertices has
$\lambda$ common neighbors, and each pair of distinct, non-adjacent
vertices has $\mu$ common neighbors.   Disjoint unions of
cliques of equal size and their complements are \emph{trivial}
examples of SR graphs.  The line graphs of complete graphs and of
complete bipartite graphs with equal parts are also SR;
we refer to them and to their complements as \emph{graphic}
SR graphs.   

The senior author has long suspected that the automorphism groups of
non-trivial, non-graphic SR graphs are ``small.''  The
following result gives a specific interpretation to this statement.

We say that a group $H$ is \emph{involved}
in a group $G$ if $G$ has subgroups\break
$L\triangleleft K\le G$
such that $H\cong K/L$ ($H$ is isomorphic to a quotient of a 
subgroup of $G$).  We note that the automorphism groups of the 
trivial and the graphic  SR graphs involve alternating groups 
of degree $\ge \sqrt{n}$.   It turns out that in all other
cases, the alternating groups involved are tiny.  This has
significant implications to attempts at subexponential
and possibly even quasi-polynomial-time isomorphism tests
for SR graphs, and it limits the possible
primitive group actions on SR graphs.  We state
the result.

\begin{theorem}[\cite{bab-srg}]  \label{thm:itcs}
If the alternating group $A_t$ is involved in the automorphism group
of a non-trivial, non-graphic SR graph with $n$ vertices then\break
 $t = O((\ln n)^2/\ln \ln n)$.
\end{theorem}

Theorem~\ref{thm:itcs} is derived from a lemma that limits the number
of fixed points of any non-identity automorphism of a regular graph in
terms of its combinatorial and spectral parameters; from this, a bound
on the orders of automorphisms follows via the Lemma above, and a
bound on $t$ is then immediate.   To apply this general result
to SR graphs, bounds on the second largest eigenvalue
and on the parameters $\lambda, \mu$ are derived from
known results.

\section{\'Akos Seress \\ (November 24, 1958 -- February 13, 2013)}
Apart from minor updates, this note was written in August 2013,
about six months after the untimely passing of \'Akos
Seress~\cite{gap}.  
I (Babai) had known \'Akos, 8 years my junior,
from his undergraduate years in the late 70s at 
E\"otv\"os University, Budapest,
where I was teaching at the time and he was a star student,
already active in research.   But our real meeting of minds
occurred at a conference in Szeged, Hungary, in summer 1986,
when by serendipity both of us missed the boat for a scenic 
afternoon ride.  It was there, on the banks of the Tisza river,  
that I began to introduce \'Akos, then a fresh
Ph.\,D. in combinatorics, to asymptotic group
theory and algorithmic group theory.
That conversation evolved into a lifelong
collaboration that produced 15 joint papers spanning a
quarter century.  I count several of our joint papers among the
best of my career; this includes our last paper~\cite{bbs}.
\'Akos became a leader in algorithmic and
computational group theory, an author of the definitive
monograph on the subject~\cite{seress-book}, a major
contributor to the GAP symbolic algebra package, 
and a speaker at ICM 2006 in the algebra section.

\'Akos was a most generous friend and colleague.  He
died at the age of 54 of renal cell carcinoma, a particularly
aggressive form of cancer, diagnosed only six months earlier.
The disease struck at the height of his creative powers, shortly 
after he had finished two breakthrough
papers, one mathematical and one computational: the above-mentioned
paper~\cite{helfgott} 
with Helfgott, a \emph{tour de force} in the combinatorial 
theory of permutation groups, and a computational work on the 
Monster group~\cite{seress-monster} that received the ``Distinguished 
paper award'' at the ISSAC 2012 conference and was hailed as 
``a groundbreaking work'' that ``marks a turning point 
in Majorana Theory.''

The Lemma discussed in this note was the fruit of the first 
hours of our collaboration, conceived even before the 
return of the boat, so I find it most appropriate to
list \'Akos as a coauthor. 

\bigskip\noindent
{\bf Acknowledgment.}  I wish to thank Ian Wanless for 
the inspiring conversation at the CameronFest last July
and for his helpful comments on earlier versions of this note.

\end{document}